\documentstyle[12pt,amscd,amssymb]{amsart}
\newtheorem{thm}{Theorem}

\newtheorem{prop}{Proposition}[section]
\newtheorem{rem}{Remark}[section]
\newtheorem{cor}{Corollary}[section]
\begin{document}

\newcommand{\C}{{\Bbb C}}
\newcommand{\HH}{{\Bbb H}}
\newcommand{\M}{{\cal M}}
\newcommand{\R}{{\Bbb R}}
\newcommand{\Z}{{\Bbb Z}}
\newcommand{\E}{{\cal E}}
\renewcommand{\O}{{\cal O}}
\renewcommand{\P}{{\Bbb P}}
\newcommand{\Q}{{\cal Q}}
\newcommand{\rk}{\o{rank}}
\title
{Moduli of Vector Bundles on Curves in Positive Characteristic}
\author{Kirti Joshi \and Eugene Z. Xia}
\date{\today} 
\address{ 
Department of Mathematics,
University of Arizona, Tucson, AZ 85721}
\email{kirti@@math.arizona.edu {\it (Joshi)}, exia@@math.arizona.edu {\it (Xia)}}
\maketitle
\begin{abstract}
Let $X$ be a projective curve of genus 2
over an algebraically closed field of characteristic 2.
The Frobenius map on $X$ induces a rational map
on the moduli space of rank-2 bundles.  We show that up
to isomorphism, there is only one (up to tensoring
by an order two line bundle) semi-stable vector 
bundle of rank 2 with determinant equal to a theta characteristic
whose Frobenius pull-back is not stable.  The indeterminacy of
the Frobenius map at this point can be resolved by introducing
Higgs bundles.
\end{abstract}

\section{Introduction and Results}
Let $X$ be a smooth projective curve of genus 2 over an algebraically
closed field $k$ of characteristic $p > 0$.  
Let $\Omega$ be its canonical bundle.
Define the (absolute)
Frobenius morphism \cite{Ka1, Ka2}
\begin{equation*}
F : X \longrightarrow X
\end{equation*}
which maps local sections $f \in \O_X$ to $f^p$.
As $X$ is smooth, $F$ is a (finite) flat map.

Let $U(r,L)$ denote the moduli of
S-equivalence classes of semi-stable vector bundles of rank
$r$ and determinant $L$ on $X$ \cite{Se}.  
Denote $J^0, J^1$ the spaces of isomorphism 
classes of line bundles of degree $0$ and $1$, respectively.
Choose a theta characteristic $L_{\theta} \in J^1$.
We study
the Frobenius pull-backs of the bundle in $S_O = U(2,\O_X)$ and 
$S_{\theta} = U(2,L_{\theta})$.
The geometry of the space $S_{\theta}$ has been studied extensively
by Bhosle \cite{Bh}.  

In general, the operation of Frobenius pull-back has a tendency to
destabilize bundles (for an example of such behavior see
\cite{Ra}).  In particular, the natural map $V \longmapsto
F^*(V)$ gives a rational map on the moduli space.  Note that the
Frobenius does give a morphism on the stack of vector bundles.

A non-semi-stable rank 2 vector bundle with trivial
determinant on $X$ must be a direct sum of line bundles with one
exception (see Proposition~\ref{prop:3.2a}).  
In Proposition~\ref{prop:3.5}, we give a necessarily and sufficient
criterion for $F^*(V)$ to be non-semi-stable in terms of theta
divisors.  If one assumes that $X$
is an ordinary curve, then the number of semi-stable vector bundles
which Frobenius destabilizes is finite (see
Proposition~\ref{prop:3.2b}).  For brevity, we formulate the following
definition.  We say a vector bundle $V$ is {\em Frobenius semi-stable
(stable)} if $V$ and $F^*(V)$ are semi-stable (stable).


For a given vector bundle $V$ on $X$, let
$$
J_2(V) = \{V \otimes L : L \in J^0, L^2 = O_X\}
$$
\begin{thm} \label{thm:2}
Suppose $p=2$.
Then there exists a bundle $V_1 \in Ext^1(L_{\theta}, \O_X)$ such
that if
$$
V \in S_{\theta} \setminus J_2(V_1),
$$
then $F^*(V)$ is semi-stable.
\end{thm}
By Theorem~\ref{thm:2}, if $V \in S_{\theta} \setminus J_2(V_1)$, 
then $F^*(V) \otimes \Omega^{-1} \in S_O$.
Hence
\begin{cor} \label{cor:1.1}
If $p=2$, then the Frobenius map induces a map
$$
F^* \otimes \Omega^{-1} : S_{\theta} \setminus J_2(V_1) \longrightarrow S_O.
$$
\end{cor}
Observe that in characteristic zero, there is no relation between 
$S_O$ and  $S_{\theta}$.  

We show that there is a natural way of resolving the indeterminacy
of the Frobenius map at the point $V_1 \in S_{\theta}$, by replacing 
$S_O, S_{\theta}$ with moduli spaces of suitable Higgs pairs.
Let $U(2,\O_X, L)$ and $U(2,L_{\theta}, L)$ be the moduli of 
semi-stable Higgs bundles with a fixed line bundle $L$ \cite{Ni}.  
There are similar notions of stability and semi-stability associated 
with a Higgs bundle \cite{Hi, Ni}.  A Higgs bundle
$(V, \phi)$ with line bundle $L$ consists of a bundle $V$ and 
a morphism $\phi : V \longrightarrow V \otimes L$.
In general, the smaller is the degree of $L$, the smaller is the
moduli $U(2,L_{\theta}, L)$.
Denote by $S_O(\Omega)$ and $S_{\theta}(L_{\theta})$ the spaces
$U(2,\O_X, \Omega)$ and $U(2,L_{\theta}, L_{\theta})$, respectively.

If $(V,\phi) \in S_{\theta}(L_{\theta})$, then 
$(F^*(V) \otimes L_{\theta}^{-1}, F^*(\phi) \otimes 1) \in S_O(\Omega)$,
where 1 denotes the constant automorphism of $L_{\theta}^{-1}$.

We say a Higgs bundle $(V, \phi)$ is {\em Frobenius semi-stable (stable)}
if $(V, \phi)$ and $(F^*(V), F^*(\phi))$ are semi-stable (stable).
\begin{thm} \label{thm:3}
Suppose $p=2$.
\begin{enumerate}
\item If $(V, \phi) \in S_{\theta}(L_{\theta})$, then either $V \in S_{\theta}$
or $V \in J_2(O_X \oplus L_{\theta})$.
\item There exist Higgs fields $\phi_0$ and $\phi_1$
such that
$(W, \phi_0)$ and $(V, \phi_1)$
are Frobenius semi-stable for all
$W \in J_2(O_X \oplus L_{\theta})$ and $V \in J_2(V_1)$.
\end{enumerate}
\end{thm}

\begin{cor}
The Frobenius defines a map on
a Zariski open set $U \subset S_{\theta}(L_{\theta})$ 
$$
F^* \otimes \Omega^{-1} : U \longrightarrow S_O(\O_X),
$$
where $U$ contains the set $S_{\theta} \setminus J_2(V_1)$ and the points
$(W, \phi_0), (V, \phi_1)$ for any $W \in J_2(O_X \oplus L_{\theta})$ and 
$V \in J_2(V_1)$.
\end{cor}

Cartier's theorem \cite{Ka1} gives a criterion for descent under
Frobenius.  We are concerned with ascent under the Frobenius.  
Higgs bundles appear
naturally in characteristic $p > 0$ context.  To see this, let $(V, \nabla)$
be a vector bundle with a (flat) connection,
$$
\nabla : V \longrightarrow \Omega \otimes_{O_X} V.
$$
One associates to the pair $(V, \nabla)$ its $p$-curvature 
(see \cite{Ka2}) which is a homomorphism of $O_X$-modules \cite{Ka2}:
$$
\psi : V \longrightarrow F^*(\Omega) \otimes_{O_X} V.
$$
Thus the pair $(V, \nabla)$ gives a Higgs pair with line bundle
$F^*(\Omega)$.  We call the pair $(V, \psi)$ a Frobenius-Higgs pair.
Moduli space of Frobenius-Higgs pairs exists (see \cite{Ni}).
>From this point of view, we are led to the consideration
of Higgs bundles to study the properties of the Frobenius on the 
moduli spaces of vector bundles and Higgs bundles.

\centerline{\sc Acknowledgments}
We thank Professor Usha Bhosle for reading a previous 
version and for her comments and suggestions for improvement. 
We thank Professors Minhyong Kim, N. Mohan Kumar, V. B. Mehta and S. Ramanan
for insightful discussions and comments.

\section{Bundle Extensions and the Frobenius Morphism}
Suppose $V$ is a vector bundle on $X$.  The slope of
$V$ is defined as
$$
\mu(V) = \deg(V) / \mbox{rank}(V).
$$

If $L$ is an line bundle on $X$, then $F^*(L) = L^p$.
The push-forward, $F_*(O_X)$, is an $O_X$-module of rank $p$
and one has the exact sequence of $O_X$-modules
\begin{equation}
0 \longrightarrow O_X \longrightarrow F_*(O_X)
\longrightarrow B_1 \longrightarrow 0
\end{equation}
where the cokernel $B_1$ is of rank $p-1$ and slope $g-1$ (see \cite{Ra}).

Tensoring the sequence with a line bundle $L$ 
gives us 
\begin{equation}
0 \longrightarrow L \longrightarrow F_*(O_X) \otimes L 
\longrightarrow B_1 \otimes L \longrightarrow 0.
\end{equation}
By the projection formula, this is
\begin{equation*}
0 \longrightarrow L \longrightarrow F_*(L^p)
\longrightarrow B_1 \otimes L \longrightarrow 0.
\end{equation*}
The associated long cohomology sequence is
\begin{equation}
... \longrightarrow H^0(B_1 \otimes L) \longrightarrow H^1(L)
\stackrel{f_L}{\longrightarrow} H^1(F_*(L^p)) \longrightarrow ...
\end{equation}
Compute the Leray spectral sequence \cite{Gr}
\begin{equation*}
\left\{
\begin{array}{lll}
E_{\infty} & \Longrightarrow & H^*(L^p)\\[2ex]
E_2^{p,q} & = & H^p(R^qF_*(L^p))
\end{array}
\right.
\end{equation*}
Since $F$ is affine, the sequence degenerates at $E_2$.
Hence 
\begin{equation*}
H^1(F_*(L^p)) \cong H^1(L^p).
\end{equation*}
Substituting this into the long exact sequence, one obtains
\begin{equation*}
... \longrightarrow H^0(B_1 \otimes L) \longrightarrow H^1(L)
\stackrel{f_L}{\longrightarrow} H^1(L^p) \longrightarrow ...
\end{equation*}

Suppose $L_1, L_2$ are two invertible sheaves.  Then the 
space of all extensions of the form
\begin{equation*}
0 \longrightarrow L_1 \longrightarrow V
\longrightarrow L_2 \longrightarrow 0
\end{equation*}
is $Ext^1(L_2, L_1)$ and is isomorphic to $H^1(L_2^{-1} \otimes L_1)$.
\begin{rem} \label{rem:2.1}
By Serre duality, one has
$$
Ext^1(L_2, L_1) = H^1(L_2^{-1} \otimes L_1) = 
H^0(\Omega \otimes L_2 \otimes L_1^{-1}).
$$
\end{rem}

Suppose $V \in Ext^1(L_2, L_1)$, i.e. $V$ fits into an exact sequence
\begin{equation*}
0 \longrightarrow L_1 \longrightarrow V
\longrightarrow L_2 \longrightarrow 0.
\end{equation*}
Since $F$ is a flat morphism, we have
\begin{equation*}
0 \longrightarrow F^*(L_1) \longrightarrow F^*(V)
\longrightarrow F^*(L_2) \longrightarrow 0.
\end{equation*}
This gives a map
\begin{equation*}
F^* : Ext^1(L_2, L_1) \longrightarrow Ext^1(F^*(L_2), F^*(L_1)) \cong
Ext^1(L_2^p, L_1^p).
\end{equation*}

\begin{prop} \label{prop:2.1}
The diagram
\begin{equation*}
\begin{CD}
H^1(L_2^{-1} \otimes L_1)  @>f_L>>       H^1(L_2^{-p} \otimes L_1^p)\\
@VV{id}V                                 @VV{id}V\\
Ext^1(L_2, L_1)            @>F^*>>     Ext^1(L_2^p, L_1^p)  
\end{CD}
\end{equation*}
commutes.
\end{prop}
\begin{pf}
This is a direct consequence of the fact that the functors $\Gamma(X,.)$
and $Hom(O_X,.)$ are equivalent.
\end{pf}

\begin{cor} \label{cor:2.1}
$F^*(V) = L_1^p \oplus L_2^p$ if and only if $V \in \ker(F^*)$.
\end{cor}
\begin{pf}
Suppose $V \in Ext^1(L_2, L_1) = H^1(L_2^{-1} \otimes L_1)$.  The 
Corollary follows directly from the
long exact sequence
$$
... \longrightarrow H^0(B_1 \otimes L_2^{-1} \otimes L_1)
\longrightarrow H^1(L_2^{-1} \otimes L_1)
\stackrel{F^*}{\longrightarrow} H^1(L_2^{-p} \otimes L_1^p) \longrightarrow ...
$$
\end{pf}

\section{The moduli of Semi-Stable Vector and Higgs Bundles}
A rank 2 vector bundle $V$ of $O_X$-modules is called stable 
(semi-stable)
if $W \subset V$ is a proper $O_X$-submodule implies $\mu(W) < \mu(V)$ 
($\mu(W) \le \mu(V)$).  The sets $S_O$ and $S_{\theta}$ are defined to
be the sets of all $S$-equivalence classes \cite{Se} of rank 2 
locally free semi-stable 
sheaves with determinant equal to 
$O_X$ and $L_{\theta}$, respectively.  Set
$$
K = \{V \in S_O : V \mbox{ is semi-stable but not stable. }\}.
$$

Given a line bundle $L$
on $X$, a Higgs bundle $(V,\phi)$ consists of a vector bundle $V$
and a Higgs field which is an $O_X$-module morphism:
$$
\phi : V \longrightarrow V \otimes L.
$$
A Higgs field can also be considered as an element in $H^0(End(V) \otimes L)$
and the Frobenius pull-back of $\phi$ is a section
$$
F^*(\phi) \in H^0(End(F^*(V)) \otimes F^*(L)) = H^0(End(F^*(V)) \otimes L^p)
$$

A Higgs bundle $(V,\phi)$ is said to be stable (semi-stable)
if $W \subset V$ with $\phi(W) \subset W \otimes L$ being a 
proper subbundle implies $\mu(W) < \mu(V)$ 
($\mu(W) \le \mu(V)$).
The set $S_O(\Omega)$ and $S_{\theta}(L_{\theta})$ are defined to
be all $S$-equivalence classes \cite{Ni} of rank 2 semi-stable Higgs 
bundles on $X$ with $L$ being 
$\Omega$ and $L_{\theta}$, respectively.

%
Suppose $V \in K$.
Then there exists $L \in J^0$ such that
$$
0 \longrightarrow L^{-1} \longrightarrow V \longrightarrow L \longrightarrow 0.
$$
The pull-back of $V$ by Frobenius then fits into the following
sequence
$$
0 \longrightarrow F^*(L^{-1}) \stackrel{f_1}{\longrightarrow} F^*(V)
\stackrel{f_2}{\longrightarrow} F^*(L) \longrightarrow 0,
$$
which is
$$
0 \longrightarrow L^{\otimes (-p)} \stackrel{f_1}{\longrightarrow} F^*(V)
\stackrel{f_2}{\longrightarrow} L^{\otimes p} \longrightarrow 0.
$$
\begin{prop} \label{prop:3.1}
$$
F^* : K \longrightarrow K.
$$
is a well-defined morphism.
\end{prop}
\begin{pf}
Let $H \subset V$ be a subbundle of maximum degree.
If $f_2 |_H = 0$, then $H = L^{\otimes (-p)}$ and
$\deg(H) = \deg(L^{\otimes (-p)}) = 0$.  If $f_2 |_H \neq 0$,
then $\deg(H) \le \deg(L^{\otimes p}) = 0$.
\end{pf}

In general, $F^*(V)$ may not be semi-stable.
For example,
$B_1$ as well as $F^*(B_1)$ are of rank $p-1$.  
By a theorem of Raynaud, the bundle $B_1$ is always semi-stable while
$F^*(B_1)$ is not semi-stable if $p > 2$ \cite{Ra}.  

If $V \in S_O$ is stable.  Then there exists $L \in J^1$ such
that 
$$
0 \longrightarrow L^{-1} \longrightarrow V \longrightarrow L \longrightarrow 0.
$$
The pull-back of $V$ by Frobenius then fits into the following
sequence
$$
0 \longrightarrow F^*(L^{-1}) \longrightarrow F^*(V) 
\longrightarrow F^*(L) \longrightarrow 0,
$$
which is
$$
0 \longrightarrow L^{\otimes (-p)} \longrightarrow F^*(V) 
\longrightarrow L^{\otimes p} \longrightarrow 0.
$$

The spaces $Ext^1(L_{\theta}, L_{\theta}^{-1})$ and $Ext^1(L_{\theta}, \O_X)$
are, respectively, isomorphic to $H^1(L_{\theta})$ and $H^1(\Omega)$, 
hence
by Riemann-Roch, are both one dimensional.  Therefore, up to isomorphism,
there exist unique nontrivial extensions $V_0, V_1$:
$$
0 \longrightarrow L_{\theta}^{-1} \longrightarrow V_0 \longrightarrow 
L_{\theta} \longrightarrow 0,
$$
$$
0 \longrightarrow O_X \longrightarrow V_1 \longrightarrow 
L_{\theta} \longrightarrow 0.
$$

\begin{prop} \label{prop:3.2}
Suppose $L \in J^1$ and 
$V$ is an extension of $L^{-1}$ by $L$.
Then either $V = L \oplus L^{-1}$ or $V \in J_2(V_0)$
\end{prop}
\begin{pf}
By Remark~\ref{rem:2.1},
the space of extension is isomorphic to
$H^0(\Omega \otimes L^{-2})$.
Since $\deg(L) = 1$, $h^0(\Omega \otimes L^{-2})$ is equal to 1 if
$L^2 = \Omega$ and 0 otherwise.  In the first case, $L \in J_2(L_{\theta})$
and $H^1(L^{-2})$ consists of the trivial extension
and a family of isomorphic bundles.
\end{pf}

\begin{prop}\label{prop:3.2a}
If $V \in K$, then $F^*(V)$ is semi-stable.
For any $V \in S_O$, if $F^*(V)$ is not semi-stable, then $F^*(V)$ is either
a direct sum of two invertible sheaves, or $V \in J_2(V_0)$.
\end{prop}
\begin{pf}
Suppose $W=F^*(V)$ is not semi-stable.  Then there exists
line bundle $M \subset V$ such that 
$$
\deg(M) = \mu(M) > \mu(W) = 0
$$
and 
$$
0 \longrightarrow M \longrightarrow W \longrightarrow M^{-1} 
\longrightarrow 0.
$$
This implies that
$$
W \in Ext^1(M^{-1}, M) \cong H^0(\Omega \otimes M^{-2}).
$$
By the above Proposition~\ref{prop:3.2}, $W \in J_2(V_0)$ or $M \oplus M^{-1}$.
\end{pf}

A curve $X$ is called ordinary if the induced map
$$
F^* : H^1(\O_X) \longrightarrow H^1(F^*(\O_X)) = H^1(\O_X)
$$
is an isomorphism.
The following proposition as a consequence of the 
results in \cite{Me, Na} was communicated to Joshi by V.B. Mehta:
\begin{prop}\label{prop:3.2b}
If $X$ is ordinary and $p > 2$, then there exists a finite set S,
such that $F^*(V)$ is semi-stable for all $V \in S_O \setminus S$.
In other words, $F^*$ induces a morphism:
$$
F^* : S_O \setminus S \longrightarrow S_O.
$$
\end{prop}
\begin{pf}
By a theorem of Narasimhan-Ramanan, when $p=0$,
$S_O \cong {\mathbb P}^3$.  Moreover, as was remarked to one of
us by Ramanan, the proof given in \cite{Na} works
in all characteristic $p \neq 2$ \cite{Na}.
If $X$ is ordinary, then by a theorem of Mehta-Subramanian \cite{Me}, 
the Frobenius morphism is \'{e}tale on a non-empty open set $U$:
$$
F^* : U \longrightarrow S_O.
$$
The open set $U$ contains $K$ which is
an effective ample divisor in $\P^3$.  
Therefore $S_O \setminus U$ is of co-dimension 3, 
hence, is a finite set.  Note that $K$ can also be
identified with the Kummer surface of $J^0$ in $\P^3$ \cite{Na}.
\end{pf}

Although unable to identify explicitly this finite set
upon which the Frobenius is not defined, we provide the 
following criterion.

\begin{prop} \label{prop:3.5}
Suppose $V \in S_O$ Then $F^*(V)$ is
not semi-stable if and only if $F^*(V)$ is an extension
$$
0 \longrightarrow M \longrightarrow F^*(V) \longrightarrow M^{-1} 
\longrightarrow 0,
$$
where $M \in J_2(L_{\theta})$.
\end{prop}

\begin{pf}
If $F^*(V)$ is indeed such an extension then $F^*(V)$ is clearly not
semi-stable as degree of $L_\theta$ is positive. We use inseparable
descent to prove the other direction. Suppose $F^*(V)$ is not
semi-stable.  Then we have an exact sequence
$$ 
0 \longrightarrow M \longrightarrow F^*(V) \longrightarrow M^{-1} \to 0,
$$
where $\deg(M) > 0$. 

Following \cite{Ka1}, $F^*(V)$ has a
connection. The the second fundamental form of the
connection is a morphism
\begin{equation}
T_X\to {\rm Hom}(M,M^{-1})=M^{-2}.
\end{equation}
As $V$ is not semi-stable this
morphism must not be the zero morphism.  Such a morphism provides 
non-zero map (after dualizing)
\begin{equation}
M^2 \longrightarrow \Omega,
\end{equation}
which is a non-zero section of the line bundle 
$\Omega \otimes M^{-2}$.  This implies
$\deg(\Omega \otimes M^{-2}) \ge 0$.  Since
$\deg(M) > 0$ and $\deg(\Omega) = 2$, we must have $\deg(M) = 1$.
Hence
$\Omega \otimes M^{-2}$ is a line bundle of degree zero which
has a section if and only if it is trivial.
Hence $\Omega = M^2$ and so $M \in J_2(L_{\theta})$.
\end{pf}

\section{The Moduli Spaces in Characteristic 2}
In this section, we follow the notations and definitions
introduced in Section 3 and assume $p=2$.
In characteristic $2$, $B_1$ is invertible and equal to 
a theta characteristic \cite{Ra}, and we choose our $L_{\theta}$ 
to be $B_1$.

\subsection{The moduli of semi-stable sheaves}
By a theorem in \cite{Na}, if $V \in S_{\theta}$, then
there exists
$L_1 \in J^0, L_2 \in J^1$ and $V$ is an extension of $L_2$ by $L_1$:
$$
0 \longrightarrow L_1 \longrightarrow V
\longrightarrow L_2 \longrightarrow 0
$$
with $L_1 \otimes L_2 = L_{\theta}$.
In other words,
$$
V \in Ext^1(L_2, L_1) = H^1(L_2^{-1} \otimes L_1).
$$
Since $L_{\theta} = B_1$, $h^0(B_1 \otimes L_2^{-1} \otimes L_1)$ is
1 if $L_{\theta} = L_2 \otimes L_1^{-1}$ and $0$ otherwise.  

Suppose
$L \in J_2(O_X)$, 
$L_1 = L$ and $L_2 = L \otimes L_{\theta}$.
By 
Corollary~\ref{cor:2.1}, 
$$
\dim(\ker(F^*)) = 1
$$ 
implying there is a unique (up to scalar) $0 \neq V \in 
Ext^1(L_2, L_1)$ such that $F^*(V) \in  O_X \oplus \Omega$.  
In fact, such an extension $V$ is isomorphic to $L \otimes V_1$.

Suppose $V \neq L \otimes V_1$.
Then Corollary~\ref{cor:2.1} implies that
$F^*(V) \neq F^*(L_1) \oplus F^*(L_2) = L_1^2 \oplus L_2^2$.
If $M \subset F^*(V)$ is a destabilizing subbundle, i.e. $\deg(M) \ge 2$, 
then $M^{-1} \otimes L_2^2$ has a global section implying 
$\deg(M) \le \deg(L_2^2) = 2$.  Moreover if $\deg(M) = 2$,
then $M = L_2^2$ implying that $F^*(V)$ contains 
$L_2^2$ as a subbundle.
Then the sequence
$$
0 \longrightarrow L_1^2 \otimes \Omega \longrightarrow F^*(V)
\longrightarrow L_2^2 \longrightarrow 0,
$$
splits.  This is a contradiction;
hence, we conclude that $\deg(M) < 2$ and $F^*(V)$ is semi-stable.
This proves Theorem~\ref{thm:2}.

\subsection{Restoring Frobenius Stability: Higgs Bundles}
Suppose $V \in S_{\theta}$.  Then the Higgs pair $(V, 0) \in S_{\theta}(L_{\theta})$.  
Hence $S_{\theta} \subset S_{\theta}(L_{\theta})$ by the map $V \longmapsto (V,0)$.  
If $(V, \phi) \in S_{\theta}(L_{\theta})$
and $V$ is not semi-stable, then $V$ is an extension
$$
0 \longrightarrow L_1 \longrightarrow V 
\stackrel{f}{\longrightarrow} L_2 \longrightarrow 0,
$$ 
where $\deg(L_1) \ge 1 > \deg(L_2)$.  Moreover $\phi(L_1)$
is not contained 
in $L_1 \otimes L_{\theta}$
(otherwise $\phi(L_1) \subset L_1 \otimes L_{\theta}$
implying $(V,\phi)$ is not semi-stable).
This implies that
there exists line bundle $L_1 \neq H \subset V$ such that 
$\phi(L_1) \in H \otimes L_{\theta}$.  Hence 
$$
\deg(L_1) \le \deg(H) + \deg(L_{\theta}).
$$
Since $L_1 \neq H$, $0 \neq f(H) \subset L_2$ implies that 
$\deg(H) \le \deg(L_2)$.  To summarize, we have the following inequalities:
$$
\deg(L_2) + \deg(L_{\theta}) \ge \deg(H) + \deg(L_{\theta}) \ge \deg(L_1) > \deg(L_2).
$$
Since $\deg(L_{\theta}) = 1$, $\deg(L_1) = \deg(H) + 1 = \deg(L_2) + 1 = 1$.  
The degree of $H$ is thus zero implying that $f(H) = L_2$, 
so the exact sequence splits.  
In addition, since $0 \neq \phi(L_1) \in H \otimes L_{\theta}$, $\phi |_{L_1}$
must be a non-zero constant morphism and
$$
L_1 = L_2 \otimes L_{\theta}.
$$
Since $L_1 \otimes L_2 = L_{\theta}$, it must be the case that
$V \in J_2(O_X \oplus L_{\theta})$.  This proves the first part of 
Theorem~\ref{thm:3}.

Given a Higgs bundle $(V, \phi) \in S_{\theta}(L_{\theta})$ with 
$\phi \in H^0(End(V) \otimes L_{\theta})$, the Frobenius pull-back 
$F^*(\phi)$ is in $H^0(End(F^*(V)) \otimes \Omega)$.  If 
$V \in S_{\theta} \setminus J_2(V_1)$, then $F^*(V)$ is semi-stable implying 
$(F^*(V), F^*(\phi))$ is semi-stable.
Suppose $L \in J_2(O_X)$.

\noindent{\sl Case 1:} $W = L \oplus L \otimes L_{\theta}$.  

We take
the Higgs field $\phi_0$ to be the identity map:
$$
1 = \phi_0 : L \otimes L_{\theta} \longrightarrow L \otimes L_{\theta}.
$$
If $M \subset L \oplus L \otimes L_{\theta}$,
then either $M = L \otimes L_{\theta}$ or $\mu(M) < 
\mu(L \oplus L \otimes L_{\theta})$.
Since $L \otimes L_{\theta}$ is not $\phi_0$-invariant, 
$(W,\phi_0)$ is stable.
The Frobenius pull-back 
$$
F^*(W) = O_X \oplus \Omega
$$
and $F^*(\phi_0)$ is again a constant map
$$
F^*(\phi_0) : \Omega \longrightarrow O_X \otimes \Omega.
$$
Again if $M \subset O_X \oplus \Omega$ then either 
$M = \Omega$ or $\mu(M) < \mu(O_X \oplus \Omega)$.
Since $\Omega$ is not $F^*(\phi_0)$-invariant, 
$(F^*(W), F^*(\phi_0))$ is stable.

\noindent{\sl Case 2:} $V = L \otimes V_1$.

The bundle $V$ is a non-trivial bundle extension of $L \otimes L_{\theta}$ 
by $L$
$$
0 \longrightarrow L \stackrel{f_1}{\longrightarrow} V 
\stackrel{f_2}{\longrightarrow} L \otimes L_{\theta} \longrightarrow 0.
$$
Tensor the sequence with $L_{\theta}$ gives
$$
0 \longrightarrow L \otimes L_{\theta} \stackrel{g_1}{\longrightarrow} V
\otimes L_{\theta} 
\stackrel{g_2}{\longrightarrow} L \otimes \Omega \longrightarrow 0.
$$
Set 
$$
\phi_1 = g_1 \circ \phi_0 \circ f_2.
$$
The Frobenius pull-back decomposes $V_1$:
$$
F^*(V) = F^*(W) = O_X \oplus \Omega.
$$
Pulling back the exact sequences by Frobenius gives
$$
0 \longrightarrow O_X \stackrel{F^*(f_1)}{\longrightarrow} O_X \oplus \Omega 
\stackrel{F^*(f_2)}{\longrightarrow} \Omega \longrightarrow 0
$$
$$
0 \longrightarrow O_X \otimes \Omega \stackrel{F^*(g_1)}{\longrightarrow} 
(O_X \otimes \Omega) \otimes \Omega
\stackrel{F^*(g_2)}{\longrightarrow} \Omega \otimes \Omega \longrightarrow 0
$$
One must show that $(V,\phi_1)$ is Frobenius stable.
If $M \subset O_X \oplus \Omega$ then either $M = \Omega$ or
$\mu(M) < \mu(O_X \oplus \Omega)$.  Hence to demonstrate 
the Frobenius stability of $(V,\phi_1)$, one only has
to show that $\Omega$ is not $F^*(\phi_1)$-invariant.

The Frobenius pull-back of $\phi_1$ is a composition:
$$
F^*(\phi_1) = F^*(g_1) \circ F^*(\phi_0) \circ F^*(f_2).
$$
Since the map $F^*(f_2)$ is surjective, the restriction map
$F^*(f_2)|_{\Omega}$ is an isomorphism.
The map $\phi_0$ is an isomorphism and $g_1$ is injective; hence,
$g_1 \circ \phi_0$ is injective.  This implies $F^*(g_1) \circ F^*(\phi_0)$
is injective.  Therefore $F^*(\phi_1)|_{\Omega}$ is injective.  
Since $\deg(\Omega) < \deg(\Omega \otimes \Omega)$, $F^*(\phi_1)|_{\Omega}$
being injective implies
$$
F^*(\phi_1)(\Omega) \not\subset \Omega \otimes \Omega \subset 
(O_X \oplus \Omega) \otimes \Omega.
$$
In other words,
$\Omega \subset O_X \oplus \Omega$ is not $F^*(\phi_1)$-invariant.
Hence $(F^*(V), F^*(\phi_1))$ is stable.  This proves Theorem~\ref{thm:3}.

\end{document}